\documentclass[12pt,platino,a4,amstex]{article}
\usepackage{amssymb}
\usepackage{epsfig}
\usepackage{array}
\usepackage{epsfig}
\usepackage[reqno]{amsmath}
\usepackage{amsthm}
\usepackage{graphicx}

\numberwithin{equation}{section}

\newtheorem{thm}{Theorem}[section]
\newtheorem{lem}{Lemma}[section]

\newcommand{\eps}{\varepsilon}
\begin{document}

\begin{center}
\Large\bf {Parameter-uniform fitted mesh higher order finite difference
scheme for singularly perturbed problem with an interior turning
point}
 \vskip .3cm
Vikas Gupta,\footnote{Department of Mathematics, The LNM Institute of Information Technology Jaipur, 302031 (India), vikasg.iitk@gmail.com, vikasg@lnmiit.ac.in} 
Sanjay K. Sahoo\footnote{Department of Mathematics, The LNM Institute of Information Technology Jaipur, 302031 (India), 16pmt003@lnmiit.ac.in}
and Ritesh K. Dubey\footnote{Research Institute, SRM University, Chennai (India), riteshkd@gmail.com}
\end{center}

\begin{abstract} In this paper, a parameter-uniform fitted mesh finite difference scheme is constructed and analyzed for a class of singularly perturbed interior turning point problems. The solution of this class of turning point problem possess two outflow exponential
boundary layers. Parameter-explicit theoretical bounds on the
derivatives of the analytical solution are given, which are used
in the error analysis of the proposed scheme. The problem is
discretized by a hybrid finite difference scheme comprises of
midpoint-upwind and central difference operator on an appropriate
piecewise-uniform fitted mesh. An error analysis has been carried out
for the proposed scheme by splitting the solution into regular and
singular components and the method has been shown second order uniform
convergent except for a logarithmic factor with respect to the
singular perturbation parameter. Some relevant numerical examples
are also illustrated to verify computationally the theoretical
aspects. Numerical experiments show that the proposed method gives competitive results in comparison to 
those of other methods exist in the literature.\newline 
\end{abstract}

\section{Introduction}\label{intro} Singularly perturbed problems arise often in the modeling of various modern complicated processes, such as viscous flow problems with large Reynolds
numbers~\cite{hirsch}, convective heat transport problems with
large P\'{e}clet numbers~\cite{jacob}, drift diffusion equation of
semiconductor device modelling~\cite{polak}, electromagnetic field
problems in moving media~\cite{hahn}, financial
modelling~\cite{black} and turbulence models~\cite{launder} etc.
Most of the singularly perturbed problems cannot be completely solved by analytical techniques. Consequently, numerical
techniques are getting much attention to get some useful insights
on the solutions of singularly perturbed problems. In general, two
classes of methods, namely, {\emph{ fitted operator methods}} and
{\emph{ fitted mesh methods}} have been used to solve such
problems.

\par Those singularly perturbed convection-diffusion problems, in which the convection coefficient vanishes at some points of the domain of the problem, are called singularly perturbed turning
point problems (SPTPPs), and zeros of the convection coefficient
are said to be \emph{turning points}. Here, we consider the
following class of singularly perturbed   two-point boundary value
problems with an interior turning point at
$x=0$~\cite{kadalbajoo:2010}:
  \begin{equation}\label{eq:intro1}
  \begin{cases}
  \mathcal{L}u(x)\equiv \varepsilon u''(x)+ a(x)u'(x) - b(x)u(x)= f(x),
  \qquad x\in \Omega = (-1, 1),\\
  u(-1)=A, \qquad u(1)=B,
  \end{cases}
  \end{equation}
  where $\varepsilon$ is a small perturbation parameter satisfying  $0< \varepsilon
  <<1$, $A$ and $B$ are given constants, $a, b$ and $f$ are sufficiently smooth functions. We impose the following restriction to ensure that the solution of Eq.~(\ref{eq:intro1})
  exhibits twin boundary layers
  \begin{equation}\label{eq:intro2}
   a(0)=0, \qquad a'(0)< 0.
  \end{equation}
Moreover, for some constant $\delta >0$ there exists a positive
constant $\alpha$, such that
\begin{equation}\label{eq:intro3}
|a(x)|\geq \alpha >0, \qquad  \delta \leq|x|\leq 1.
\end{equation}
Also $b(x)$ is required to be bounded below by some positive
constant $\beta$, $i.e.$,
\begin{equation}\label{eq:intro5}
b(x)\geq \beta>0, \qquad x\in \bar\Omega =[-1,1],
\end{equation}
to guarantee that the operator $\mathcal{L}$ is inverse monotone on $[-1,1]$ and to exclude the so-called resonance phenomena \cite{ackerberg:1970}. We also impose the following restriction to
ensure that there are no other turning points in the interval
$[-1, 1]$:
\begin{equation}
|a'(x)| \geq  \left|\frac{a'(0)}{2}\right|, \qquad x\in \bar\Omega
= [-1, 1].\label{eq:intro6}
\end{equation}
This class of singularly perturbed turning point problem (SPTPP)
~(\ref{eq:intro1}) has a unique solution possess twin outflow
boundary layers of exponential type at both end points $x=\pm 1$,
under the assumptions ~(\ref{eq:intro2})-~(\ref{eq:intro6}).

\par It is very difficult to deal singularly perturbed turning point problems analytically. The study of these problems received
much attention in the literature due to the complexity involved in
finding uniformly valid asymptotic expansions unlike non-turning
problems. Some authors, such as, Jingde~\cite{jingde:1991},
O'Malley~\cite{malley:1974, malley:1991}, Wasow~\cite{wasow}
studied qualitative aspects of these problems, namely, existence,
uniqueness and asymptotic behavior of the solution.

\par In general, since the convection coefficient has zero inside the domain therefore numerical treatment of singularly perturbed
turning point problem becomes more difficult than the singularly perturbed non-turning point problems. Abrahamsson~\cite{abrahamsson:1977},
Berger et al.~\cite{berger:1984} and Farrell~\cite{farrell:1988}
establish {\emph{a priori}} bounds for interior turning point
problems; in particular it is shown that a bound is independent of
singular perturbation parameter $\varepsilon$ if and only if
reaction coefficient is greater than zero at the turning point. It
is also shown there how the ratio of reaction coefficient
$b(x)\geq 0$ and first derivative of convection coefficient
$a'(x)$, $i.e.$, $\lambda =b(x)/a'(x)$ at the turning point plays
a key role in determining the behavior of the
solution~\cite{berger:1984}. It is shown that for $\lambda <0$, the
solution is smooth near turning point and two outflow boundary
layers of exponential type exhibits at both the endpoints of the
domain. In this case the turning point is sometimes called a
{\emph{diverging flow}} or {\emph{expansion}} turning point. On
the other hand, if $\lambda>0$, there is in general no boundary
layers exhibited and an {\emph{interior layer}} appears at the
turning point, the nature of which depends in a fundamental way on
$\lambda$. For $0<\lambda <1$, the interior layer is called
{\emph{cusp layer}} because it can be approximately modelled by a
cusp-like function. Interior layer turning point is sometimes
called a {\emph{converging flow}} or {\emph{compression}} turning
point. In inviscid fluid dynamics, the
diverging flow turning point corresponds to a {\emph{sonic point}}
and converging flow turning point to a {\emph{shock point}}.  For
the case $b(x)=0$ at the turning point, the solution exhibits a
very interesting phenomenon called {\emph{Ackerberg-O'Malley
resonance phenomenon}}~\cite{ackerberg:1970}.

\par Berger et al.~\cite{berger:1984} also show that the modified version of El Mistikawy Werle scheme
is uniformly convergent of  $O(h^{\min(\lambda, 1)})$ in the
$L^{\infty}[-1,\,1]$ norm using the analytic bounds obtained
in~\cite{berger:1984}. Farrell~\cite{farrell:1988} obtained a set
of sufficient conditions for uniform convergence in the discrete
$L^{\infty}$ norm on uniform mesh, not only for exponentially
fitted schemes, but also for a large class of schemes of upwinded
type. Kadalbajoo and Patidar~\cite{kadalbajoo:2001} gave a
numerical scheme based on cubic spline approximation with
nonuniform mesh for SPTPP~(\ref{eq:intro1})-(\ref{eq:intro6}) and
established second order $\varepsilon$-uniform convergence.
Natesan et al.~\cite{natesan:2003} proposed a numerical method
based on the classical upwind finite difference scheme on a
Shishkin mesh and proved that the proposed scheme is uniformly
convergent of almost order one. In~\cite{kadalbajoo:2010},
Kadalbajoo and Gupta derived asymptotic bounds for the derivatives
of the analytical solution of
SPTPP~(\ref{eq:intro1})-(\ref{eq:intro6}) and proposed a
computational method comprises B-spline collocation scheme on a
non-uniform Shiskin mesh. They shown that this scheme is second
order accurate in the maximum norm. Kadalbajoo et
al.~\cite{kadalbajoo:2011} also suggested B-spline collocation
with artificial viscosity on uniform mesh for the same class of
SPTPP~(\ref{eq:intro1})-(\ref{eq:intro6}). In~\cite{munyakazi:2014}, Munyakazi and Patidar conclude that convergence acceleration Richardson extrapolation technique on existing numerical schemes for the above class of turning point problem does not improve the rate of convergence. However, Becher and Roos~\cite{becher:2015} show that Richardson extrapolation on upwind scheme with piecewise-uniform Shishkin mesh works fine and improves the accuracy to $O(N^{-2}\ln^2 N)$ under the assumption $\varepsilon \leq CN^{-1}$. Recently, Munyakazi et al.~\cite{munyakazi:2019} proposed a fitted operator finite difference scheme for singularly perturbed turning point problem having an interior layer and also shown that with Richardson extrapolation technique, accuracy and order of convergence of the scheme can be improved upto two.  For a general review of existing literature on asymptotic and numerical analysis of turning point problems, one can see~\cite{sharma:2013}.
\par
In this paper, we focus to devise a second order uniformly
convergent finite difference scheme for SPTPP~(\ref{eq:intro1}) on
piecewise uniform mesh of Shishkin type without using any convergence acceleration technique like Richardson extrapolation. The proposed method
combines the midpoint upwind difference scheme and classical
central finite difference scheme on piecewise uniform mesh. The
requirements of higher order truncation error and monotonicity
play a vital role in the construction of this scheme. One can
observe the fact that the classical central difference scheme is
monotone if $\varepsilon$ is relatively large than the convection coefficient $a$ $i.e.$, if
$\varepsilon \geq Ch||a||$, where $h$ is the mesh width and has second order truncation error
on uniform mesh. On the other hand, midpoint upwind difference
operator
 is monotone for all value of $\varepsilon$ and for relatively large convection coefficient $a$ than the reaction coefficient $b$ such
 that $h||b||\leq C\alpha$. Moreover, midpoint upwind operator possess second order truncation error away from the boundary layer region.
 Also, Shiskin mesh equally distribute the number of mesh points inside and outside the boundary layers, therefore one can gets a coarse mesh region outside the boundary layer
 and fine mesh region inside the boundary layer. Utilizing these facts, we employ midpoint upwind difference scheme in coarse
mesh region and central difference operator in fine mesh region of
Shishkin mesh. Since, central difference operator yields first
order truncation error at transition points, we use midpoint
upwind operator on transition points. Such type of higher order
scheme for singularly perturbed non-turning convection-diffusion
problem was introduced by Stynes and Roos~\cite{stynes:1997}. To
analyze the proposed scheme theoretically, we split the numerical
solution into regular and singular components and analyze them
separately by using tools such as
truncation error bounds, discrete
minimum principle and appropriate choices of barrier functions.\\
{\bf Notation.} Throughout the paper we use $C$ as a generic
positive constant independent of $\varepsilon$ and mesh
parameters.  For any given function $g(x)\in C^{k}(\bar\Omega)$
($k$ a non-negative integer), $||g||$ is a global maximum norm
over the domain $\bar\Omega$ defined by
\begin{equation*}
||g||=\max_{\bar\Omega}|g(x)|.
\end{equation*}
\section{A-\emph{priori} Estimates for Continuous Problem}
In this section some bounds of the exact solution and its
derivatives are discussed. These bounds will be needed for error
analysis of proposed numerical scheme in later sections.
Derivation of these bounds are well known and can be found
in~\cite{kadalbajoo:2010}. Systematically, we use minimum
principle to derive these bounds.
\begin{lem}
\label{lem:maxprin}(\cite{kadalbajoo:2010}.) {\bf (Minimum
Principle)} Let $\psi(x) \in C^{2}(\bar\Omega)$ and $\psi(\pm
1)\geq 0$. Then $\mathcal{L}\psi(x)\leq 0, \forall x\in \Omega$
implies that $\psi(x)\geq 0, \forall x\in \bar\Omega.$
\end{lem}
Since the concerned SPTPP~(\ref{eq:intro1})-(\ref{eq:intro6}) is
linear, minimum principle ensure the existence and uniqueness of
the classical solution. Using the above minimum principle, one can
easily prove the following uniform stability estimate for the
differential operator $\mathcal{L}$.
\begin{lem}\label{stability estimate} (\cite{kadalbajoo:2010}.) {\bf (Uniform Stability Estimate)}  $\forall \varepsilon
>0 $, solution $u(x)$ of the SPTPP~(\ref{eq:intro1})-(\ref{eq:intro6}), satisfies the
following stability estimate:
\begin{equation*}
||u(x)||\leq \frac{||f||}{\beta} +\max(|A|,|B|), \qquad \forall
x\in \bar\Omega.
\end{equation*}
\end{lem}
To exclude the turning point $x=0$ and to obtain the bounds for
the solution $u$ and its derivatives in the non-turning point
region of the domain, we divide the domain $\bar\Omega$ into three
subdomains as $\Omega_{1}=[-1,-\delta]$,
$\Omega_{2}=[-\delta,\delta]$ and $\Omega_{3}=[\delta,1]$ such
that $\bar \Omega =\Omega_{1}\cup \Omega_{2} \cup \Omega_{3}$,
where $0<\delta \leq 1/2$. Further, following theorem gives bounds
for the derivatives of $u$ in the subintervals $\Omega_1$ and
$\Omega_3$ individually.
\begin{thm} (\cite{kadalbajoo:2010}.)\label{thm:contin1}
If $a,\, b$ and $f$ $\in C^{m}(\bar \Omega), m>0$, then solution
$u(x)$ of the SPTPP~(\ref{eq:intro1})-(\ref{eq:intro2}) satisfies
the following bounds for any $\delta >0$:
\begin{equation*}
|u^{j}(x)| \leq C\left(1+\varepsilon
^{-j}\exp\left(-\frac{\alpha(1+x)}{\varepsilon}\right)\right),
\quad j=1,\cdots ,m+1, \quad x\in \Omega_{1},
\end{equation*}
\begin{equation*}
|u^{j}(x)| \leq C\left(1+\varepsilon
^{-j}\exp\left(-\frac{\alpha(1-x)}{\varepsilon} \right)\right),
\quad j=1,\cdots ,m+1, \quad x\in \Omega_{3},
\end{equation*}
\end{thm}
Next, we state a theorem, which gives the bounds for the
derivatives of the solution in the turning point region $\Omega_2$
and deduce that the solution is smooth in subdomain $\Omega_2$.
\begin{thm}\label{thm:contin2}(\cite{berger:1984}.)
Let $u(x)$ be the solution of SPTPP defined from the
equations(\ref{eq:intro1})-(\ref{eq:intro6}), and $a,\, b, \,f\in
C^{m}(\bar\Omega), m>0$. Then for $\varepsilon >0$ and
sufficiently small $\delta >0$, there exists a positive constant
$C$ such that
\begin{equation*}
|u^{(j)}(x)|\leq C, \quad j=1,2,\ldots, m, \quad \forall x\in
\Omega_{2}.
\end{equation*}
\end{thm}
It turns out that the bounds for continuous solution $u(x)$ given
in Theorem~\ref{thm:contin1} and Theorem~\ref{thm:contin2} are not
adequate to obtain $\varepsilon$-uniform error estimate for the
proposed scheme. Therefore, to analyze the proposed scheme
correctly, we need to derive more precise bounds on these
derivatives by decomposing the solution into regular component $v$
and singular component $w$ as
\begin{equation*}
u(x)=v(x)+w(x), \qquad  \forall x\in \bar\Omega,
\end{equation*}
where the smooth component $v$ satisfies homogeneous problem
$\mathcal{L}v(x)=f(x)$ and singular component satisfies
homogeneous problem $\mathcal{L}w(x)=0$ with appropriate boundary
conditions. Using the technique given in~\cite{kadalbajoo:2010}, we get the following bounds for smooth and singular components in the region $\Omega_1$:
\begin{equation*}
|v^{(j)}(x)|\leq C(1+\varepsilon^{((m-1)-j)}
e^{-\alpha(1+x)/\varepsilon}),  \qquad \forall x\in \Omega_{1},
\end{equation*}
\begin{equation*}
|w^{(j)}(x)|\leq C\varepsilon^{-j}
e^{-\alpha(1+x)/\varepsilon}.\qquad \qquad \qquad \forall x\in
\Omega_{1}.
\end{equation*}
In the same manner, we can obtain analogous
estimates for subinterval $\Omega_{3}$, while the solution $u(x)$
and its derivatives are smooth in the subinterval $\Omega_{2}$. Hence, on the whole domain $\bar\Omega$, the bounds on $v$ and $w$, and their
derivatives are given in the following theorem:
\begin{thm}\label{thm:contin3}(\cite{kadalbajoo:2010}.)
Let $a, b$ and $f\in C^{m}(\bar \Omega), m>0,$  then for all $j,~
0\leq j\leq m,$ the smooth component satisfies
\begin{equation*}
|v^{(j)}(x)|\leq
C\left(1+\varepsilon^{((m-1)-j)}\left(\exp\left(-\frac{\alpha(1+x)}{\varepsilon}\right)+\exp\left(-\frac{\alpha(1-x)}{\varepsilon}
\right)\right)\right), \quad \forall x\in \bar\Omega,
\end{equation*}
and the singular component satisfies
\begin{equation*}
|w^{(i)}(x)|\leq
C\varepsilon^{-i}\left(\exp\left(-\frac{\alpha(1+x)}{\varepsilon}\right)+\exp\left(-\frac{\alpha
(1-x)}{\varepsilon}\right)\right), \quad \forall x\in \bar\Omega.
\end{equation*}
\end{thm}
\section{Fitted Mesh Higher-Order Scheme}\label{discrete}
In this section, first we construct fitted piecewise-uniform mesh
$\bar\Omega^N$ of Shishkin type to discretize the domain
$\bar\Omega$ and then employ a specially designed finite
difference scheme on this mesh to discretize the
SPTPP~(\ref{eq:intro1})-(\ref{eq:intro2}). The fitted mesh
$\bar\Omega^N$  is constructed by dividing $\bar\Omega$ into three
subintervals $\Omega_L=[-1, -1+\tau], \Omega_{C}=[-1+\tau,
1-\tau]$ and $\Omega_{R}=[1-\tau, 1]$ such that $\bar\Omega=
\Omega_{L} \cup \Omega_{C} \cup \Omega_{R}$. For $N\geq
2^r,\,r\geq 3$ be an integer, $\bar\Omega^N$ divides each of the
subintervals $\Omega_L$ and $\Omega_R$ into $N/4$ mesh intervals
and $\Omega_C$ with $N/2$ mesh intervals such that
$\bar\Omega^{N}=\{-1=x_0,\,x_1,\ldots,x_{N/4}=-1+\tau,\ldots
x_{3N/4}=1-\tau,\ldots ,x_N=1\}$. Here, the transition parameter
is obtained by taking
\begin{equation*}
\tau =\min \left\lbrace \frac{1}{4}, \tau_0\varepsilon \ln N\right
\rbrace.
\end{equation*}
The constant $\tau_0$ is independent of the parameter
$\varepsilon$ and the number of mesh points $N$ and will be chosen
later on during the analysis of proposed scheme. This mesh is
coarse on $\Omega_C$ and fine on $\Omega_L$ and on $\Omega_R$. If
$h$ and $H$ are fine and coarse mesh width respectively, then mesh
width $h_i=x_i-x_{i-1},\,i=1,\ldots,N,$ is defined as
\begin{equation*}\label{eq:mesh1}
h_i=
\begin{cases}
h=4\tau/N, &\text{ $i=1,2,\ldots, N/4$},\\
H=4(1-\tau)/N, &\text{ $i=N/4+1,\ldots,3N/4$},\\
h=4\tau/N, &\text{ $i=3N/4+1,\ldots, N$.}
\end{cases}
\end{equation*}
One can easily observe that
\begin{equation*}
N^{-1}\leq H\leq  4 N^{-1},\quad h=4\tau_0 \varepsilon N^{-1}\ln
N<N^{-1},\quad  H+h=4N^{-1}.
\end{equation*}
Since, convection coefficient $a(x)$ changes its sign at the
turning point $x=0$, therefore, we construct a finite difference
scheme $\mathcal{L}^NU_i=\tilde{f_i},\,i=1,2,\ldots,N-1$ to
discretize the SPTPP~(\ref{eq:intro1}) in the following manner
\begin{equation}\label{eq:hybrid1}
\mathcal{L}^{N}U (x_i) \equiv
\begin{cases}
\mathcal{L}_{\text{c}}^N U\equiv\varepsilon\delta^{2}U_i+a_{i}D^{0}U_i-b_{i}U_i=f_{i}, & i=1,2,\ldots, N/4-1,\\
\mathcal{L}_{\text{mp}}^N U \equiv \varepsilon\delta^{2}U_i+a_{i\pm 1/2}D^{\pm}U_i-(bU)_{i\pm 1/2}=f_{i\pm 1/2}, & i=N/4,\ldots,3N/4,\\
\mathcal{L}_{\text{c}}^N U \equiv\varepsilon\delta^{2}U_i+a_{i}D^{0}U_i-b_{i}U_i=f_{i}, & i=3N/4+1,\ldots, N-1,\\
U_0=A,\qquad U_{N}=B,
\end{cases}
\end{equation}
where,
\begin{equation*}
\mathcal{L}_{\text{mp}}^N U\equiv
\begin{cases}
\varepsilon\delta^{2}U_i+a_{i+ 1/2}D^{+}U_i-(bU)_{i+
1/2}=f_{i+ 1/2}, &\text{if}~a_i>0\\
\varepsilon\delta^{2}U_i+a_{i- 1/2}D^{-}U_i-(bU)_{i- 1/2}=f_{i-
1/2},&\text{if}~a_i<0.
\end{cases}
\end{equation*}
Here, we used the following definition to construct above scheme
\begin{equation*}
v_i= v(x_i),\quad v_{i+1/2}=\frac{v_i+v_{i+1}}{2},\quad
v_{i-1/2}=\frac{v_{i-1}+v_i}{2}, \quad
\widehat{h}_i=\frac{h_i+h_{i+1}}{2},
\end{equation*}
\begin{equation*}
D^{+}v_{i}=\frac{v_{i+1}-v_{i}}{h_{i+1}},~,D^{-}v_{i}=\frac{v_{i}-v_{i-1}}{h_{i}},~D^{0}v_{i}=\frac{v_{i+1}-v_{i-1}}{2\widehat{h}_{i}},
~\delta^{2}v_{i}=\frac{\left(D^{+}v_{i}-D^{-}v_{i}\right)}{
\widehat{h}_{i}}.
\end{equation*}
It is clear that proposed finite difference operator
$\mathcal{L}^N$ in scheme~(\ref{eq:hybrid1}) is a combination of
central difference operator $\mathcal{L}^N_{\text{c}}$ and
midpoint upwind difference operator $\mathcal{L}^N_{\text{mp}}$,
which is constructed
 by using knowledge judiciously about the sign of the convection term, location of the turning point and truncation error behavior of these operators. After simplifying the terms
in~(\ref{eq:hybrid1}), the difference scheme takes the form
$\mathcal{L}^NU_i\equiv
p_i^{l}U_{i-1}+p_i^{c}U_i+p_i^{r}U_{i+1}=\tilde{f_i}$, where the
coefficients are given by
\begin{align*}
&p_{i}^{l}=\left(\frac{\varepsilon}{h_{i}\widehat{ h}_{i}}-\frac{a_{i}}{2\widehat{h}_{i}}\right),\qquad \qquad p_{i}^{c}=\left(-p_i^{l}-p_i^{r}-b_i\right),\\
&p_{i}^{r}=\left(\frac{\varepsilon}{h_{i+1}\widehat{ h}_{i}}+\frac{a_{i}}{2\widehat{h}_{i}}\right),~i=1,2,\ldots,N/4-1,~3N/4+1\ldots,N-1,\\
&p_{i}^{l}=\left(\frac{\varepsilon}{h_{i}\widehat{ h}_{i}}\right),~p_{i}^{c}=\left(-p_i^{l}-p_i^{r}-b_{i+1/2}\right),~p_{i}^{r}=\left(\frac{\varepsilon}{h_{i+1}\widehat{ h}_{i}}+\frac{a_{i+1/2}}{h_{i+1}}-\frac{b_{i+1}}{2}\right),~\text{if}~a_i>0\\
&p_{i}^{l}=\left(\frac{\varepsilon}{h_{i}\widehat{
h}_{i}}-\frac{a_{i-1/2}}{h_{i}}-\frac{b_{i-1}}{2}\right),~p_{i}^{c}=\left(-p_i^{l}-p_i^{r}-b_{i-1/2}\right),~p_{i}^{r}=\left(\frac{\varepsilon}{h_{i+1}\widehat{
h}_{i}}\right),~\text{if}~a_i<0,\\
& i=N/4,\ldots,3N/4.
\end{align*}
\section{Uniform Convergence}
Here, in this section first we shall establish the consistency and
stability estimate through discrete minimum principle and then
analyze proposed numerical method~(\ref{eq:hybrid1}) for
$\varepsilon$-uniform convergence by analogous decomposition of
discrete solution into smooth and singular components as of
continuous solution.
\begin{lem}\label{lem:uniform1}{\bf( Discrete Minimum Principle)}
Let us suppose that $N\geq N_0$, where
\begin{equation}\label{eq:uniform1}
\frac{h||a||}{2\varepsilon} < 1,
~~i.e.,~~2\tau_0||a||<\frac{N_0}{\ln
N_0},~~and~~\frac{2||b||}{N_0}\leq \alpha.
\end{equation}
Then the operator $\mathcal{L}^N$ defined by~(\ref{eq:hybrid1})
satisfies a discrete minimum principle, $i.e.,$ if $\psi^N$ is a
mesh function that satisfies $\psi^N_0\geq 0,\,\psi^N_N\geq 0$ and
$\mathcal{L^N}\psi_i^N\leq 0,$ for $1\leq i\leq N-1$, then
$\psi_i^N\geq 0$ for $0\leq i\leq N$.
\end{lem}
{\bf Proof.} In order to establish the discrete minimum principle, We simply check that the associated system matrix is $M$-matrix with the choice of the midpoint upwind and central
difference operator used in the definition of the difference
scheme~(\ref{eq:hybrid1}). It allow us to establish the following
inequalities on the coefficients of the difference operator
$\mathcal{L}^N$:
\begin{equation}\label{eq:uniform2}
p_i^l>0,\qquad p_i^r>0,\qquad p_i^l+p_i^c+p_i^r<0,\qquad
i=1,2,\ldots,N-1.
\end{equation}
 In the case
of central difference operator $\mathcal{L}^N_c$, conditions in
(\ref{eq:uniform2}) are satisfied if $h||a||<2\varepsilon$ $i.e.,$
if $N_0(\ln N_0)^{-1}>2\tau_0||a||$, then one can check $p_i^l>0$
and $p_i^r>0,\,p_i^l+p_i^c+p_i^r<0$ for $1\leq i\leq N/4-1$ and
$3N/4+1\leq i\leq N-1$. For the case of midpoint upwind operator
$\mathcal{L}^{N}_{\text{mp}}$, the conditions in
(\ref{eq:uniform2}) are satisfied if $H||b||<2\alpha$ $i.e.,$ if
$2||b||<\alpha N_0.$ From these sign patterns on the coefficients
of associated system matrix, one can deduce that operator
$\mathcal{L}^N$ is of negative type and therefore satisfies a
discrete minimum principle. Moreover, it ensures that the operator
is uniformly stable in the maximum norm.
\begin{lem}\label{lem:uniform2}
Let $Z^N_i$ be any mesh function such that $Z^N_0=Z^N_N=0.$ Then
for all $i,\,0\leq i\leq N,$ we have
\[|Z^N_i|\leq \frac{1}{\beta}\max_{1\leq j\leq N-1}|\mathcal{L}^NZ^N_j|.\]
\end{lem}
{\bf Proof.} Let us introduce two comparison functions defined by
\begin{equation*}
\Psi^{\pm}_i= \frac{1}{\beta}\max_{1\leq j\leq
N-1}|\mathcal{L}^NZ^N_j|\pm Z^N_i.
\end{equation*}
Clearly one can notice that $\Psi^N_0=\Psi^N_N \geq 0,$ since
$Z^N_0=Z^N_N=0.$ Furthermore, for $1\leq i\leq N-1,$ we have
\begin{equation*}
\mathcal{L}^N \Psi^{\pm}_i= -\frac{b}{\beta}\max_{1\leq j\leq
N-1}|\mathcal{L}^NZ^N_j|\pm \mathcal{L}^NZ_i^N \leq 0,
\end{equation*}
as $b/\beta \geq 1.$ Therefore, discrete minimum principle
(\ref{lem:uniform1}) implies that $\Psi^{\pm}_i\geq 0,\,0\leq
i\leq N$, which gives desired result.
\par Further, using the valid Taylor's series expansion, we obtained the following truncation error estimates for different finite difference operator employed in the operator $\mathcal{L}^N$: On a uniform mesh with step size $\tilde{h}$, we have
\begin{equation*}
|\mathcal{L}_{\text{c}}^Nu_i-(\mathcal{L}u)(x_i)|\leq C(\varepsilon \tilde{h}^2|u^{(iv)}|+\tilde{h}^2|u^{(iii)}|).
\end{equation*}
On an arbitrary non-uniform mesh, we have
\begin{equation*}
|\mathcal{L}_{\text{c}}^Nu_i-(\mathcal{L}u)(x_i)|\leq C\left(\varepsilon (h_i+h_{i+1}) |u^{(iii)}|+(h_i+h_{i+1})|u^{(ii)}|\right).
\end{equation*}
Here, one can notice that order of truncation error is  reduced to one only if the central difference operator is employed on arbitrary non-uniform mesh
instead of uniform mesh. Moreover, We have the following truncation error bounds corresponding to the midpoint upwind difference operator, which are valid for both uniform and non-uniform mesh:
\begin{equation*}
|\mathcal{L}_{\text{mp}}^Nu_i-(\mathcal{L}u)(x_{i-1/2})|\leq
\begin{cases}
C\left(\varepsilon (h_i+h_{i+1})|u^{(iii)}|+h_{i+1}^2(|u^{(iii)}|+|u^{(ii)}|+|u^i|)\right), & \text{if}~a(x)>0,\\
C\left(\varepsilon (h_i+h_{i+1})|u^{(iii)}|+h_{i}^2(|u^{(iii)}|+|u^{(ii)}|+|u^i|)\right), & \text{if}~a(x)<0.
\end{cases}
\end{equation*}
Note that the order of truncation error is higher by one in the convection term for midpoint upwind operator than the centered difference operator on a non-uniform mesh. This is the reason to apply midpoint upwind scheme at the transition points $(-1+\tau)$ and $(1-\tau)$ of proposed mesh.
\par Further the solution $U$ of the discrete problem can be decomposed in an analogous manner as that of
the continuous solution $u$ into the following sum
\begin{subequations}\label{eq:uniform3}
\begin{equation}
U=V+W,
\end{equation}
where,
\begin{equation}
\mathcal{L}^NV=f,\quad V(-1)=v(-1),\,V(1)=v(1),
\end{equation}
\begin{equation}
\mathcal{L}^NW=0,\quad W(-1)=w(-1),\,W(1)=w(1).
\end{equation}
\end{subequations}
Therefore, the error can be written in the form
\begin{equation*}
U - u = (V-v) + (W-w),
\end{equation*}
so the errors in the smooth and singular components of the solution can be estimated separately.
\begin{lem}\label{lem:uniform3}{\bf (Error in smooth component)}
Assume that $N\geq N_0$ satsifies the assumption~(\ref{eq:uniform1}). Then the regular component of the error satisfies the following error bound
\begin{equation*}
|(V-v)(x_i)|\leq
\begin{cases}
CN^{-2},\quad &\forall i=0,1,\ldots,N/4-1,3N/4+1,\ldots, N,\\
CN^{-1}(\varepsilon+N^{-1}) &\forall i=N/4,N/4+1,\ldots,3N/4.
\end{cases}
\end{equation*}
\end{lem}
{\bf Proof.} Using the usual truncation error estimates given above and bounds for the smooth component $v$ given in Theorem~(\ref{thm:contin3}), we have
\begin{align*}
|\mathcal{L}^N(V-v)(x_i)|&\leq
\begin{cases}
CN^{-2}(\varepsilon |v^{(iv)}|+|v^{(iii)}|), \qquad \qquad \forall i=0,1,\ldots,N/4-1,3N/4+1,\ldots, N,\\
CN^{-1}(\varepsilon |v^{(iii)}|+N^{-1}(|v^{(iii)}|+|v^{(ii)}|+|v^i|)),\,\forall i=N/4,N/4+1,\ldots,3N/4.\\
\end{cases}\\
&\leq \begin{cases}
CN^{-2},\quad &\forall i=0,1,\ldots,N/4-1,3N/4+1,\ldots, N,\\
CN^{-1}(\varepsilon+N^{-1}), &\forall i=N/4,N/4+1,\ldots,3N/4,
\end{cases}
\end{align*}
and applying Lemma~\ref{lem:uniform2}, we obtain the required result.
\par Since $a_i\geq \alpha>0,\,\forall x_i<0,\,i=1,\ldots,N/2$ and $a_i\leq -\alpha<0,\,\forall x_i>0,\,i=N/2+1,\ldots,N-1$, we consider both the region $[-1,0]$ and $[0,1]$ individually to get the error estimates for the layer component $(W-w)$. Therefore, we consider the following barrier functions for a positive constant $\gamma$:
\begin{equation}\label{eq:uniform4}
\Phi_i^L=\begin{cases}
\prod_{j=1}^i \left(1+\frac{\gamma h_j}{\varepsilon}\right)^{-1},\,&i=1,\ldots N/2,\\
1,&i=0.
\end{cases} \Phi_i^R=\begin{cases}
\prod_{j=i+1}^N \left(1+\frac{\gamma h_j}{\varepsilon}\right)^{-1},\,&i=N/2,\ldots N-1,\\
1,&i=N.
\end{cases}
\end{equation}
\par First we prove the following technical result.
\begin{lem}\label{lem:uniform4}
If $2\gamma <\alpha$, the barrier functions satsisfy the inequalities
\begin{equation*}
\mathcal{L}^N\Phi_i^L\leq 0,\quad \forall i=1,2,\ldots,N/2,\quad \mathcal{L}^N\Phi_i^R\leq 0, \quad \forall i=N/2,\ldots,N-1.
\end{equation*}
\end{lem}
{\bf Proof.} We begin with the left hand barrier function $\Phi_i^L$ and analyze each of the different discretizations used in the definition of the operator $\mathcal{L}^N$. First, in the case of midpoint upwind operator with $a(x)>0$, we have $\mathcal{L}^N\Phi_i^L=\mathcal{L}_{\text{mp}}^N\Phi_i^L=\varepsilon\delta^{2}\Phi_i^L+a_{i+ 1/2}D^{+}\Phi_i^L-(b\Phi^L)_{i+1/2}$. Using the properties, $\Phi_i^L>0,$ $\,D^{+}\Phi_i^L=-\frac{\gamma}{\varepsilon}\Phi_{i+1}^L<0,$ and $\delta^2\Phi_i^L=\left(\frac{\gamma}{\varepsilon}\right)^2\frac{h_{i+1}}{\widehat{h}_i}\Phi_{i+1}^L>0$, and with the condition  $2\gamma <\alpha$, one can easily observe that
\begin{align*}
\mathcal{L}_{\text{mp}}^N\Phi_i^L&=\left(\frac{\gamma^2}{\varepsilon}\frac{h_{i+1}}{\widehat{h}_i}-a_{i+1/2}\frac{\gamma}{\varepsilon}-\frac{b_{i+1}}{2}\right)\Phi_{i+1}^L-\frac{b_i}{2}\Phi_i^L\\
&=\left(2\frac{\gamma^2}{\varepsilon}\left(\frac{h_{i+1}}{2\widehat{h}_i}-1\right)+\left(2\frac{\gamma^2}{\varepsilon}-a_{i+1/2}\frac{\gamma}{\varepsilon}-\frac{b_{i+1}}{2}\right)-\frac{b_i}{2}\left(1+\frac{\gamma h_{i+1}}{\varepsilon}\right)\right)\Phi_{i+1}^L\leq 0.
\end{align*}
In the case of central difference operator with $a(x)>0$, we have
\begin{align*}
\mathcal{L}_{\text{c}}^N\Phi_i^L&= \left(2\frac{\gamma^2}{\varepsilon}\left(\frac{h_{i+1}}{2\widehat{h}_i}-1\right)+\left(2\frac{\gamma^2}{\varepsilon}-a_i\frac{\gamma}{\varepsilon}\frac{h_{i+1}}{2\widehat{h}_i}\right)\right)\Phi^L_{i+1}-\left(a_i\frac{\gamma}{\varepsilon}\frac{h_{i}}{2\widehat{h}_i}+b_i\right)\Phi_i^L\leq 0.
\end{align*}
Similarly, applying the midpoint upwind operator for the case $a(x)<0$, we have
\begin{equation*}
\mathcal{L}^N_{\text{mp}}\Phi_i^R=\left(2\frac{\gamma^2}{\varepsilon}\left(\frac{h_{i}}{2\widehat{h}_i}-1\right)+\left(2\frac{\gamma^2}{\varepsilon}+a_{i-1/2}\frac{\gamma}{\varepsilon}-\frac{b_{i-1}}{2}\right)-\frac{b_i}{2}\left(1+\frac{\gamma h_{i}}{\varepsilon}\right)\right)\Phi_{i-1}^L\leq 0.
\end{equation*}
In the same manner if we use central difference operator with $a(x)<0$, we also get $\mathcal{L}^N_{\text{c}}\Phi_i^R\leq 0.$ It completes the proof.
\begin{lem}\label{lem:uniform5}
The barrier functions $\Phi_i^L$ and $\Phi_i^R$ and layer component $W$ satisfy
\begin{align*}
|W_i|\leq C\Phi_i^L,\quad \forall i=0,1,\ldots N/2,\quad |W_i|\leq C\Phi^R_i, \quad \forall i=N/2,\ldots, N.
\end{align*}
Moreover, following bounds are valid for the layer component $W$ in no layer region $\Gamma_C$
\begin{equation*}
|W_i|\leq CN^{-2},\qquad \forall i=N/4,\ldots,3N/4.
\end{equation*}
\end{lem}
{\bf Proof.} Construct the barrier functions $\Psi_L^{\pm}(x_i)=C\Phi_i^L\pm W_i,\,i=0,\,1,\ldots,N/2$. By Lemma~\ref{lem:uniform4}, we have
$\mathcal{L}^N\Psi_L^{\pm}(x_i)\leq 0$. Now using the discrete minimum principle we obtain the requred bound. Furthermore, to obtain the bound for $W_i$ in no layer region $[-1+\tau, 0]$, we have for $i>N/4$:
\begin{align*}
\Phi_i^L\leq \Phi_{\frac{N}{4}}^L&=\prod_{j=1}^{N/4}\left(1+\frac{\gamma h_j}{\varepsilon}\right)^{-1}=\left(1+\frac{\gamma h}{\varepsilon}\right)^{-N/4}=
\left(1+\frac{4\gamma \tau}{\varepsilon N}\right)^{-N/4}=(1+4\gamma \tau_0N^{-1}\ln N)^{-N/4}\\
&=(1+8N^{-1}\ln N)^{-N/4}=((1+8N^{-1}\ln N)^{-N/8})^2\leq CN^{-2},
\end{align*}
for the choice of $\tau_0=2/\gamma$. Here, we have used the inequality $\ln(1+t)>t(1-t/2)$ with $t=8N^{-1}\ln N$ to prove $(1+8N^{-1}\ln N)^{-N/8}\leq 8N^{-1}$.
Using similar argument for barrier function $\Phi_i^R$, we obtain desired bounds for $W_i$ in the domain $[0,\,1]$.
\begin{lem}\label{lem:uniform6}{\bf (Error in singular component)}
Assume that $N\geq N_0$ satsifies the assumption~(\ref{eq:uniform1}) and $2\gamma <\alpha$. Then the singular component of the error satisfies the following
error estimates
\begin{equation*}
|(W-w)(x_i)|\leq
\begin{cases}
CN^{-2}(\ln N)^2,\quad &\forall i=0,1,\ldots,N/4-1,\,3N/4+1,\ldots, N,\\
CN^{-2}, &\forall i=N/4,N/4+1,\ldots,3N/4.
\end{cases}
\end{equation*}
\end{lem}
{\bf Proof.} We split our discussion into the two cases of boundary layer region $\Omega_L\cup \Omega_R$ and
no boundary layer region $\Omega_C$  to analyze the singular component of the error. Since $\Omega_C=[-1+\tau,\,0]\cup [0,\,1-\tau]$, it is sufficient to
consider only the subinterval $[-1+\tau,\,0]$ and using same
argument one can get similar estimate for the subinterval $[0,\,1-\tau]$. Both $w$ and $W$ are small in $\Omega_C$, therefore we will use
triangle inequlaity, Theorem~\ref{thm:contin3}, Lemma~\ref{lem:uniform5} instead of the usual truncaton error argument, to get the required error bounds on layer component in $[-1+\tau,\,0]$.
For $i=N/4,\ldots,N/2$, using triangle inequality, we have
\begin{align}\label{eq:uniform5}
 |(W-w)(x_i)|&\leq |W(x_i)|+|w(x_i)|\nonumber\\
 &\leq C\prod_{j=1}^{i}\left(1+\frac{\gamma h_j}{\varepsilon}\right)^{-1}+C\exp \left(-\frac{\alpha (1+x_i)}{\varepsilon}\right)\nonumber\\
 &\leq C\prod_{j=1}^{i}\left(1+\frac{\gamma h_j}{\varepsilon}\right)^{-1}\quad \text{(since $\displaystyle{e^{-\alpha(1+x_i)/\varepsilon}\leq \Phi_i^L}$)}\nonumber\\
 &\leq CN^{-2} \quad \text{(Using Lemma~\ref{lem:uniform5}).}
\end{align}
Proceeding in a similar manner in subinterval $[0,\,1-\tau]$, one can prove
\begin{equation}\label{eq:uniform6}
  |(W-w)(x_i)|\leq CN^{-2},\qquad \forall i=N/2,\ldots,3N/4.
\end{equation}
We now consider the boundary layer region $\Omega_L$ to estimate the singular component of the error. In this case, we obtain the following singular component of
 the local truncation error estimates for $i=1,2,\ldots,N/4-1$:
 \begin{align}\label{eq:uniform7}
  |\mathcal{L}^N(W-w)(x_i)|&\leq Ch^2\left(\varepsilon |w^{(iv)}|+|w^{(iii)}|\right) \nonumber\\
  &=16 CN^{-2}\tau^2\left(\varepsilon |w^{(iv)}|+|w^{(iii)}|\right)\nonumber\\
  &\leq CN^{-2}\varepsilon^2(\ln N)^2\left(\varepsilon^{-3}\exp\left(-\frac{\alpha(1+x_i)}{\varepsilon}\right)\right)\nonumber\\
  &=C\left(\frac{N^{-2}(\ln N)^2}{\varepsilon}\exp\left(-\frac{\alpha(1+x_i)}{\varepsilon}\right)\right)\nonumber\\
  &\leq C\left(\frac{N^{-2}(\ln N)^2}{\varepsilon}\Phi_i^L \right).
 \end{align}
From the Eq.~(\ref{eq:uniform5}), $|(W-w)(x_{N/4})|\leq CN^{-2}$, also we have $|(W-w)(x_{0})|=0$. Therefore, if we choose
\begin{equation*}
 \Psi^{\pm}(x_i)=CN^{-2}\left(1+(\ln N)^2\Phi_i^L\right)\pm (W-w)(x_i),\quad \forall i=0,1,\ldots,N/4.
\end{equation*}
as our barrier functions, one can easily see that both the functions satisfy $\Psi^{\pm}(x_0)\geq 0,$ and $\Psi^{\pm}(x_{N/4})\geq 0$. Moreover, $\mathcal{L}^N
\Psi^{\pm}(x_i)=-Cb_iN^{-2}+CN^{-2}(\ln N)^2\mathcal{L}^N\Phi_i^L\pm \mathcal{L}^N(W-w)(x_i)\leq 0$ by Lemma~\ref{lem:uniform4} and estimate given in Eq.~(\ref{eq:uniform7}).
Therefore, by applying discrete minimum principle, we obtain $\Psi^{\pm}(x_i)\geq 0,\,\forall i=0,1,\ldots,N/4$, which gives
\begin{equation}\label{eq:uniform8}
 |(W-w)(x_i)|\leq CN^{-2}(1+(\ln N)^2\Phi_i^L),\qquad \forall i=0,1,\ldots,N/4.
\end{equation}
Now to get the bounds for $\Phi_i^L$ for $i=1,2,\ldots,N/4$, we use the approach given in~\cite{kellog}, for that we have
\begin{align*}
\Phi_i^L&=\prod_{j=1}^i\left(1+\frac{\gamma h_j}{\varepsilon}\right)^{-1}=\left(1+\frac{\gamma h}{\varepsilon}\right)^{-i}=\left(1+\frac{\gamma h}{\varepsilon}\right)^{-(1+x_i)/h}=\left(1-\frac{\gamma h}{\gamma h+\varepsilon}\right)^{(1+x_i)/h},\\
\Rightarrow\ln \Phi_i^L&=\frac{(1+x_i)}{h}\ln \left(1-\frac{\gamma h}{\gamma h+\varepsilon}\right)\leq \frac{(1+x_i)}{h}\left(-\frac{\gamma h}{\gamma h+\varepsilon}\right)
=-\frac{\gamma(1+x_i)}{\gamma h+\varepsilon},
\end{align*}
now taking the exponential of both sides, we get the following estimates:
\begin{equation*}
 \Phi_i^L\leq \exp\left(-\frac{\gamma(1+x_i)}{\gamma h+\varepsilon}\right)=\exp\left(-\frac{
 \gamma i h}{\gamma h+\varepsilon}\right), \quad \forall i=1,2,\ldots,N/4.
\end{equation*}
Since in $\Omega_L$, we have $h_i=h=4\tau_0\varepsilon N^{-1}\ln N,\,\forall i=1,2,\ldots,N/4$, therefore from the above we lead to the following estimate
\begin{align*}
 \Phi_i^L\leq \exp\left(- \frac{4i\gamma \tau_0N^{-1}\ln N}{1+4\gamma \tau_0N^{-1}\ln N}\right)&=N^{-4iN^{-1} \gamma \tau_0/(1+4\gamma \tau_0N^{-1}\ln N)}\\
 &=N^{-8iN^{-1}/(1+8N^{-1}\ln N)},\quad \text{with the choice of $\tau_0=2/\gamma$},\\
 &=N^{-8iN^{-1}}N^{64iN^{-1}(N^{-1}\ln N)/(1+8N^{-1}\ln N)}\\
 &\leq CN^{-8iN^{-1}}, \qquad \forall i=1,2,\ldots,N/4.
\end{align*}
Thus, from the Eq.~(\ref{eq:uniform8}), we have
\begin{align*}
 |(W-w)(x_i)|&\leq CN^{-2}(1+N^{-8iN^{-1}}(\ln N)^2),\\
 &\leq C \max\left\lbrace N^{-2}, N^{-(2+8i/N)}(\ln N)^2)\right\rbrace \\
 &\leq CN^{-2}(\ln N)^2, \qquad \forall i=0,1,\ldots,N/4.
\end{align*}
Proceeding in the same manner, one can get similar estimate for singular component of the error in $\Omega_R$, $i.e.,$ for $i=3N/4+1,\ldots,N$, which completes the
proof.

\par The Lemma~\ref{lem:uniform3} and Lemma~\ref{lem:uniform6} together gives the following main result of $\varepsilon$-uniform error estimate for the proposed fited mesh finite difference scheme.
\begin{thm}\label{thm:uniform1}
 Assume that $N\geq N_0$ satsifies the assumption~(\ref{eq:uniform1}) and $2\gamma <\alpha$. Then the continuous solution $u$ of the SPTPP~(\ref{eq:intro1})-(\ref{eq:intro6})
 and discrete solution $U$ of the finite difference approximation~(\ref{eq:hybrid1}) satisfy the following $\varepsilon$-uniform error estimate:
 \begin{equation*}
  \sup_{0<\varepsilon \leq 1}||U-u||_{\bar\Omega}\leq
  \begin{cases}
  CN^{-2}(\ln N)^2,\quad &\forall i=0,1,\ldots,N/4-1,\,3N/4+1,\ldots, N,\\
  CN^{-1}(\varepsilon+N^{-1}) &\forall i=N/4,N/4+1,\ldots,3N/4.
  \end{cases}
 \end{equation*}
\end{thm}
From the above error estimates, it is clear that for $\varepsilon \leq N^{-1}$, proposed finite difference scheme is almost second order accurate upto a logarithmic factor.
\section{Numerical Results and Discussions}
\label{numerical} In this section, we apply the constructed numerical method~(\ref{eq:hybrid1}) to the following two SPTPP to demonstrate both the
accuracy and order of convergence. Both of the problems exhibit a turning point at $x=1/2$.\\
{\textbf{Example 1.}} In this test problem, we consider the
following SPTPP:
\begin{subequations}
\begin{equation}
\varepsilon u''(x)-2(2x-1)u'(x)-4u(x)=0, \qquad x\in(0,1),
\end{equation}
\begin{equation}
u(0)=1, \qquad u(1)=1.
\end{equation}\label{eq:numeric1}
\end{subequations}
The exact solution of this problem is given by
\begin{equation}\label{eq:numeric2}
u(x)=e^{-2x(1-x)/\varepsilon}.
\end{equation}
As we know the exact solution, we can exactly compute the maximum
pointwise errors for every $\varepsilon$ in the following standard
way
\begin{equation}\label{eq:numeric3}
E_{\varepsilon}^{N}=
\max_{x_{i}\in\bar\Omega_{N}}|u(x_{i})-U^{N}(x_{i})|,
\end{equation}
where superscript $N$ denotes the number of mesh points used.
Further, we compute the $\varepsilon$-uniform maximum pointwise
error using
\begin{equation}\label{eq:numeric4}
E^{N}=\max_{\varepsilon}E_{\varepsilon}^{N}.
\end{equation}
Approximation for the order of local convergence
$\rho_{\varepsilon}^N$ is obtained in the following way
\begin{equation}\label{eq:numeric5}
\rho_{\varepsilon}^N=\log_2\frac{E_{\varepsilon}^N}{E_{\varepsilon}^{2N}}.
\end{equation}
Computed numerical results and comparison with other numerical methods available in literature  are given in Tables~\ref{table1}-\ref{table2}.\\
 {\textbf{Example 2.}} This example is corresponds to the following
 nonhomogeneous SPTPP:
 \begin{subequations}\label{eq:numeric6}
 \begin{equation}
 \varepsilon u''(x)-2(2x-1)u'(x)-4u(x)=4(4x-1), \qquad x\in (0,1),
 \end{equation}
 \begin{equation}
 u(0)=1, \qquad u(1)=1.
 \end{equation}
 \end{subequations}
 Again it posses the continuous solution given by
\begin{equation}
u(x)=-2x+2e^{-2x(1-x)/\eps}+e^{-2x(1-x)/\varepsilon}erf
((2x-1)/\sqrt{2\eps})/erf (1/\sqrt{2\varepsilon}),
\end{equation}
 where the approximations for maximum pointwise errors and numerical order of convergence are estimated as
 for the  Example 1 and corresponding numerical results are displayed in
 Tables~\ref{table3}-\ref{table4}.
 \clearpage
\begin{table}[!h]
{\small \caption{Maximum pointwise errors $E^N_{\varepsilon}$ and
 order of convergence $\rho^N_{\varepsilon}$ for Example
1}\label{table1}
 \vspace*{.5cm}
\begin{tabular}{|c|c|c|c|c|c|c|c|}
\hline

$\varepsilon \downarrow$ &N=16&N=32 &N=64 &N=128 &N=256&N=512&N=1024 \\
\hline
$10^{0}$&8.9709E-3&   4.3375E-3 &  2.1245E-3&   1.0502E-3  & 5.2199E-4 &  2.6020E-4  & 1.2990E-4  \\
       &1.0484&   1.0298 &  1.0164 &  1.0086&   1.0044 &  1.0022&            \\
$10^{-1}$&1.7821E-2&   5.8482E-3 &  1.7776E-3&   9.1441E-4 &  4.6337E-4  & 2.3316E-4  & 1.1694E-4  \\
        &1.6075&  1.7180 &  0.9591 & 0.9807&   0.9908 &  0.9955   & \\
$10^{-2}$&2.6001E-2&   1.1289E-2&   4.2974E-3 &  1.5223E-3&   5.1594E-4&   1.6820E-4 &  5.3062E-5  \\
        &1.2037&   1.3934 &  1.4972 &  1.5610  & 1.6170 & 1.6644     & \\
$10^{-3}$&2.6811E-2 &  1.1489E-2  & 4.3147E-3 &  1.4985E-3 &  4.9852E-4 &  1.6066E-4 &  5.0562E-5 \\
        &1.2226&   1.4129 &  1.5258 &  1.5878 &  1.6337 &  1.6679    &    \\
$10^{-4}$&2.6891E-2 &  1.1506E-2 &  4.3123E-3&   1.4912E-3&   4.9223E-4 &  1.5649E-4 &  4.8201E-5 \\
        &1.2247 &  1.4159 &  1.5320&   1.5990&   1.6533&   1.6989      &    \\
$10^{-5}$&2.6899E-2 &   1.1508E-2 &  4.3120E-3 &  1.4904E-3 &4.9150E-4  & 1.5595E-4  & 4.7838E-5\\
        &   1.2249&   1.4162  & 1.5327 &  1.6004 &  1.6561 & 1.7049 & \\
$10^{-6}$&2.6900E-2 &  1.1508E-2 &  4.3120E-3 &  1.4903E-3 &  4.9143E-4 &  1.5590E-4  & 4.7800E-5  \\
        & 1.2249 &1.4163 &  1.5328  & 1.6005 &  1.6564 &  1.7055    &      \\
$10^{-7}$&2.6900E-2 &  1.1508E-2 &  4.3120E-3&   1.4903E-3&   4.9142E-4&   1.5589E-4 &  4.7802E-5\\
      &1.2249 &  1.4163  & 1.5328 &  1.6005 &  1.6564 &  1.7054      & \\
$10^{-8}$&2.6900E-2 &  1.1508E-2 &  4.3120E-3&   1.4903E-3 &  4.9142E-4&   1.5590E-4 &  4.7795E-5  \\
       &1.2249 &  1.4163&  1.5328&   1.6005 &  1.6563&   1.7057& \\
$10^{-9}$&2.6900E-2 &  1.1508E-2 &  4.3120E-3 &  1.4905E-3 &  4.9162E-4 &  1.5593E-4 &  4.8395E-5 \\
        &   1.2249&  1.4163 & 1.5326 &  1.6001 &  1.6566 &  1.6880
 & \\
\hline
$E_{10^{-9}}^N~\cite{natesan:2003}$&1.796E-1 &1.178E-1 &8.00E-2 &4.95E-2 &2.98E-2 &1.72E-2 &9.7E-3\\
$\rho_{10^{-9}}^N~\cite{natesan:2003}$ &0.6084 & 0.5583 &
0.6926&0.7321& 0.7929&0.8264
&\\
\hline
$E_{10^{-9}}^N~\cite{kadalbajoo:2010}$&4.7221E-2 &1.8175E-2 &6.6037E-3 &2.3400E-3 &8.2109E-4 &2.8839E-4 &9.0532E-5\\
$\rho_{10^{-9}}^N~\cite{kadalbajoo:2010}$ &1.3775 & 1.4606 &
1.4967&1.5110& 1.5096&1.6716
&\\
 \hline
\end{tabular}
}
\end{table}
\begin{table}[!h]
{\small \caption{Maximum pointwise errors $E^N_{\varepsilon}$ and
 order of convergence $\rho^N_{\varepsilon}$ for Example
1}\label{table2}
 \vspace*{.5cm}
\begin{tabular}{|c|c|c|c|c|c|c|c|}
\hline

$\varepsilon \downarrow$ &N=16&N=32 &N=64 &N=128 &N=256&N=512&N=1024 \\
\hline
$2^{-12}$&2.6879E-2 &  1.1504E-2 &  4.3127E-3 &  1.4924E-3 &  4.9336E-4 &  1.5729E-4&   4.8719E-5\\
      &1.2244 &  1.4154 &  1.5309  & 1.5970&   1.6492&   1.6909   & \\
$2^{-16}$&2.6899E-2 &  1.1508E-2  & 4.3120E-3  & 1.4904E-3 &  4.9155E-4 &  1.5599E-4  & 4.7860E-5 \\
       &1.2249&  1.4162&   1.5327&   1.6003&   1.6559&   1.7045& \\
\hline
$E_{2^{-12}}^N~\cite{kadalbajoo:2001}$&4.1E+2 & ---- &6.9E-2 &1.5E-2 &3.7E-3 &9.2E-4 &2.3E-4\\
\hline
$E_{2^{-12}}^N~\cite{kadalbajoo:2011}$&7.670E-2 & 3.465E-2&1.646E-2 &8.018E-3 &3.957E-3 &1.966E-3 &9.840E-4\\
$\rho_{2^{-12}}^N~\cite{kadalbajoo:2011}$ &1.1464 & 1.0739 &
1.0377&1.0188& 1.0091&0.9985
&\\
$E_{2^{-16}}^N~\cite{kadalbajoo:2011}$ &7.670E-2& 3.465E-2&
1.646E-2& 8.018E-3& 3.957E-3& 1.966E-3& 9.797E-4\\
$\rho_{2^{-16}}^N~\cite{kadalbajoo:2011}$& 1.1464& 1.0739& 1.0377&
1.0188 &1.0091& 1.0049&\\
 \hline
\end{tabular}
}
\end{table}
\clearpage

\begin{table}[!h]
{\small \caption{Maximum pointwise errors $E^N_{\varepsilon}$ and
order of convergence $\rho^N_{\varepsilon}$ for Example
2}\label{table3}
 \vspace*{.5cm}
\begin{tabular}{|c|c|c|c|c|c|c|c|}
\hline

$\varepsilon \downarrow$ &N=16&N=32 &N=64 &N=128 &N=256&N=512&N=1024 \\
\hline
$10^{0}$&2.3328E-2 &  1.1634E-2  & 5.8187E-3 &  2.9087E-3 &  1.4546E-3 &  7.2731E-4 &  3.6366E-4  \\
       & 1.0037  & 0.9996 &   1.0003&   0.9998&   0.9999 &  1.000&            \\
$10^{-1}$&5.4473E-2 &  1.7786E-2 &  4.9326E-3 &  2.5413E-3&   1.2883E-3 &  6.4849E-4 &  3.2531E-4  \\
        &1.6148 &  1.8503 &  0.9568&  0.9800 &  0.9904  & 0.9953 & \\
$10^{-2}$&7.8004E-2 &  3.3867E-2 &  1.2892E-2 &  4.5668E-3 &  1.5478E-3 &  5.0459E-4&   1.5919E-4 \\
        &1.2037 &  1.3934 &  1.4972 &  1.5610&   1.6170 &  1.6644     & \\
$10^{-3}$&8.0434E-2&   3.4468E-2 &  1.2944E-2&   4.4955E-3 &  1.4956E-3&   4.8197E-4&   1.5169E-4 \\
        &1.2226& 1.4129&1.5258& 1.5878&  1.6337& 1.6679 &    \\
$10^{-4}$&8.0674E-2&   3.4519E-2&   1.2937E-2&   4.4735E-3 &  1.4767E-3&  4.6947E-4&   1.4460E-4 \\
        &1.2247&   1.4159 & 1.5320 &  1.5990 &  1.6533  & 1.6989      &    \\
$10^{-5}$& 8.0698E-2 &  3.4524E-2 &  1.2936E-2 &  4.4711E-3&   1.4745E-3&   4.6786E-4 &  1.4351E-4\\
        &  1.2249 &   1.4162 &  1.5327  & 1.6004&   1.6561&   1.7049 & \\
$10^{-6}$&8.0701E-2&   3.4525E-2&   1.2936E-2&   4.4708E-3&   1.4743E-3&   4.6770E-4 &  1.4340E-4  \\
        & 1.2249 &  1.4163 &  1.5328 &  1.6005 &  1.6564 &  1.7056   &      \\
$10^{-7}$&8.0701E-2&   3.4525E-2&   1.2936E-2&   4.4708E-3&   1.4743E-3 &  4.6768E-4 &  1.4341E-4\\
      &1.2249 &  1.4163  & 1.5328 &  1.6005 &  1.6564 &  1.7054      & \\
$10^{-8}$&8.0701E-2&   3.4525E-2&   1.2936E-2 &  4.4707E-3 &  1.4740E-3 &  4.6770E-4 &  1.4304E-4  \\
       &1.2249&   1.4163 &  1.5328 &  1.6007 &  1.6561 &  1.7091& \\
$10^{-9}$&8.0701E-2&   3.4525E-2 &  1.2936E-2 &  4.4714E-3 &  1.4749E-3 &  4.6780E-4  & 1.4319E-4 \\
        &  1.2250 &  1.4163 &  1.5326 &  1.6001 &1.6566&  1.7065& \\
\hline
$E_{10^{-9}}^N~\cite{kadalbajoo:2010}$&2.4007E-1& 1.1937E-1 &5.8785E-2& 2.7630E-2 &1.1739E-2& 4.9664E-3& 1.9735E-3 \\
 $\rho_{10^{-9}}^N~\cite{kadalbajoo:2010}$& 1.0080& 1.0220 &1.0893& 1.2349& 1.2410& 1.3314& \\
 \hline
\end{tabular}
}
\end{table}
\begin{table}[!h]
{\small \caption{Maximum pointwise errors $E^N_{\varepsilon}$ and
 order of convergence $\rho^N_{\varepsilon}$ for Example
2}\label{table4}
 \vspace*{.5cm}
\begin{tabular}{|c|c|c|c|c|c|c|c|}
\hline

$\varepsilon \downarrow$ &N=16&N=32 &N=64 &N=128 &N=256&N=512&N=1024 \\
\hline
$2^{-12}$&8.0636E-2  & 3.4511E-2 &  1.2938E-2 &  4.4773E-3  & 1.4801E-3 &  4.7188E-4 &  1.4616E-4\\
      &1.2244&   1.4154&   1.5309 &  1.5970&   1.6492 &  1.6909& \\
$2^{-16}$&8.0697E-2 &  3.4524E-2 &  1.2936E-2 &  4.4712E-3 &  1.4746E-3&   4.6796E-4 &  1.4358E-4 \\
       &1.2249&   1.4162 &  1.5327 &  1.6003  &1.6559&   1.7045& \\
\hline
$E_{2^{-12}}^N~\cite{kadalbajoo:2011}$&2.557E-2 &1.155E-2& 5.485E-3& 2.673E-3& 1.319E-3& 6.553E-4& 3.280E-4\\
$\rho_{2^{-12}}^N~\cite{kadalbajoo:2011}$ &1.1466& 1.0743& 1.0370&
1.0190 &1.0092& 0.9985
&\\
$E_{2^{-16}}^N~\cite{kadalbajoo:2011}$ &2.557E-2& 1.155E-2 &5.485E-3& 2.673E-3& 1.319E-3&6.553E-4& 3.266E-4\\
$\rho_{2^{-16}}^N~\cite{kadalbajoo:2011}$& 1.1466 &1.0743& 1.0370& 1.0190 &1.0092 &1.0046&\\
 \hline
\end{tabular}
}
\end{table}
\clearpage

\begin{figure}[htp]
\begin{center}
\includegraphics[scale=0.75]{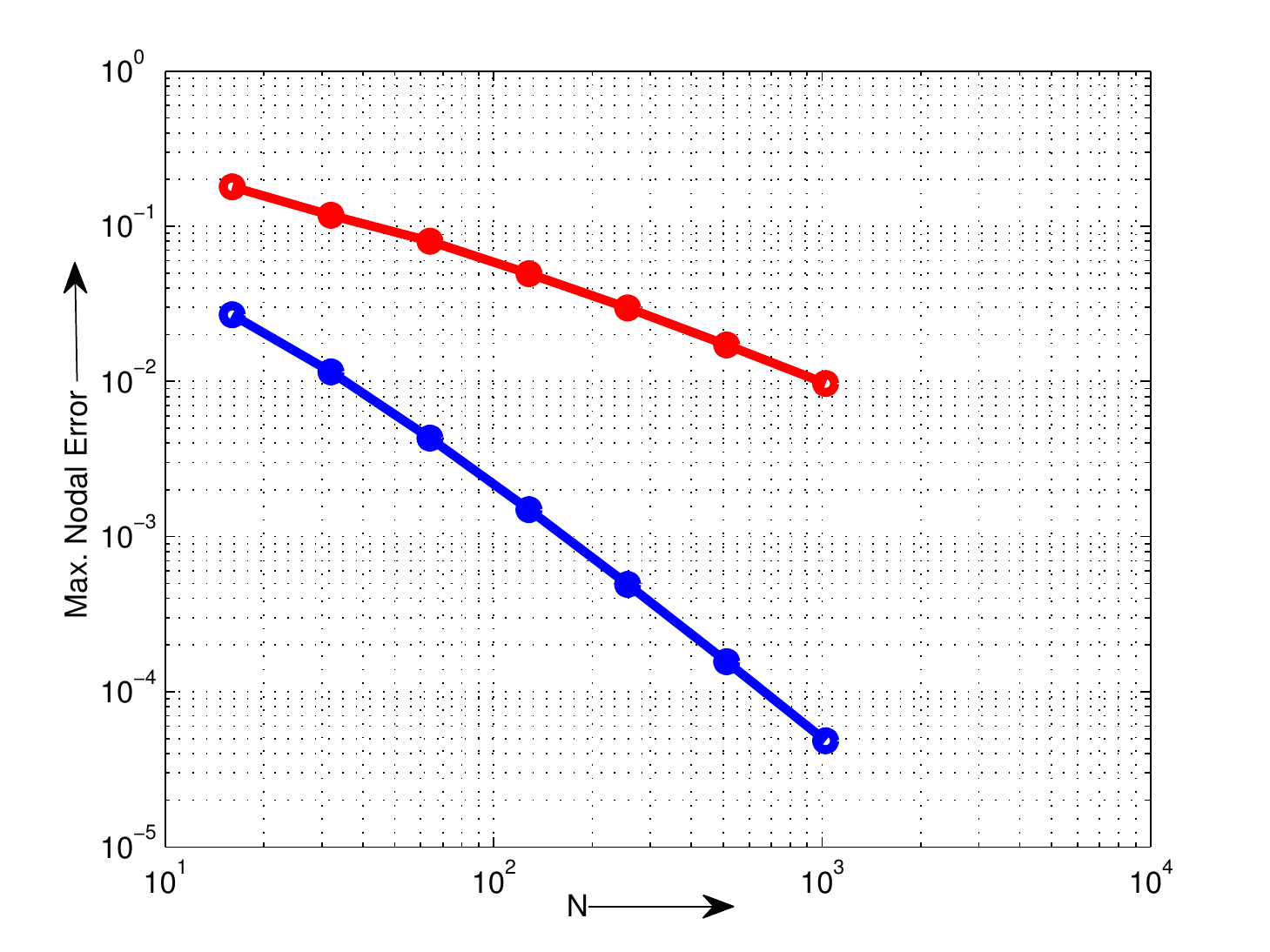}
\end{center}
\caption{{\small{Loglog plot of the maximum nodal errors with
$\varepsilon=10^{-9}$ correspond to finite difference
scheme~\ref{eq:hybrid1} (blue line) and the upwind
scheme~\cite{natesan:2003} (red line) for Example 1}}}\label{fig1}
\end{figure}

\par Numerical results presented in Tables~\ref{table1}-\ref{table4}
show that the accuracy of the proposed finite difference scheme is
in good agreement with the theoretical prediction. We apply both
the forward midpoint upwind and backward midpoint upwind operator
depending upon the sign $a(x)$ to tackle the stability of the
proposed finite difference scheme. Table~\ref{table1} and
Table~\ref{table3} display the maximum pointwise error and order
of convergence for Example 1 and Example 2 respectively for
different value of $\varepsilon$ and $N$. Table~\ref{table1} and
Table~\ref{table3} indicate that the order of convergence of
presented fitted mesh finite difference scheme~(\ref{eq:hybrid1}) is one for
$\varepsilon \geq 10^{-1}$ and almost of order two upto a
logarithmic factor for $\varepsilon <10^{-1}$. It happens because
for moderate value of $\varepsilon$, $i.e.,$ for $\varepsilon >
N^{-1}$, midpoint upwind operator is first order convergent as
given in Theorem~\ref{thm:uniform1} and in this case error
correspond to the midpoint upwind operator dominates the error
correspond to the central difference operator. Numerical results
given in Tables~\ref{table1}-\ref{table4} also show that the
maximum nodal errors decreases and order of convergence increases
as the number of mesh point increases. One can observe that as
$\varepsilon$ is getting smaller for a particular value of mesh
points $N$, both the maximum pointwise error and order of
convergence are going to stabilized.

\par A comparison given in Table~\ref{table1} for $\varepsilon=10^{-9}$, clearly
indicate that the maximum pointwise errors are much smaller and
order of convergence is much larger in this article than those
obtained in~\cite{natesan:2003} using upwind finite difference
operator. It verify numerically the theoretical estimates that
hybrid finite difference scheme~(\ref{eq:hybrid1}) is second order
$\varepsilon$-uniform convergent as opposed to the first order
uniform convergence of upwind finite difference
scheme~\cite{natesan:2003} for turning point problems. We have not
made comparison of  numerical results for Example 2 with the
finite difference scheme given in~\cite{natesan:2003} because of
authors used double mesh principle instead of analytical solution
to get pointwise errors in~\cite{natesan:2003}. Thus with almost
same computational effort, proposed finite difference scheme gives
more  accuracy and rapid convergence then the finite difference
scheme~\cite{natesan:2003}. We also compare proposed finite
difference scheme with the spline based numerical
methods~\cite{kadalbajoo:2001,kadalbajoo:2010,kadalbajoo:2011}
numerically for both the Examples~1-2 and found that present
scheme produce lesser pointwise errors and larger order of
convergence than the spline based numerical
methods~\cite{kadalbajoo:2001,kadalbajoo:2010,kadalbajoo:2011}. Furthermore, one can see, Example~1 is analogous to the Testproblem~1 
in~\cite{becher:2015} and our results are comparable to those extrapolation results in ~\cite{becher:2015} as both numerical schemes are of almost 
second order convergence $O(N^{-2}\ln (N)^2)$ under the common assumption $\varepsilon \leq CN^{-1}$ for a given number of mesh points $N$.

\par In Figure~\ref{fig1}, the Loglog graph of maximum pointwise errors is given correspond to the proposed scheme (blue line) and upwind finite difference scheme~\cite{natesan:2003}
(red line). This plot also indicate that the
error of our scheme diminishing at the rate of $1/N^2$ while error
correspond to scheme~\cite{natesan:2003} approaches to zero almost
as $1/N\rightarrow 0$. Thus all numerical evidences support our
theoretical estimates.

\end{document}